# THE OVERHAND SHUFFLE MIXES IN $\Theta(N^2 \log N)$ STEPS

By Johan Jonasson[1]

*Chalmers University of Technology*

The overhand shuffle is one of the "real" card shuffling methods in the sense that some people actually use it to mix a deck of cards. A mathematical model was constructed and analyzed by Pemantle [*J. Theoret. Probab.* **2** (1989) 37–49] who showed that the mixing time with respect to variation distance is at least of order $n^2$ and at most of order $n^2 \log n$. In this paper we use an extension of a lemma of Wilson [*Ann. Appl. Probab.* **14** (2004) 274–325] to establish a lower bound of order $n^2 \log n$, thereby showing that $n^2 \log n$ is indeed the correct order of the mixing time. It is our hope that the extension of Wilson's lemma will prove useful also in other situations; it is demonstrated how it may be used to give a simplified proof of the $\Theta(n^3 \log n)$ lower bound of Wilson [*Electron. Comm. Probab.* **8** (2003) 77–85] for the Rudvalis shuffle.

**1. Introduction.** The distribution of any aperiodic irreducible Markov chain approaches the unique stationary distribution as the number of steps increases, but how many steps are needed to get sufficiently close? This question has attracted much attention over the last two decades or so. One important reason for this is that computer development has allowed the development of powerful Markov chain Monte Carlo techniques for simulation.

Particular interest has been paid to the case when the state space is the symmetric group in which case one often refers to the Markov chain as a card shuffling chain. A great variety of different card shuffling chains have been studied. Among these are many simplified versions of "real" shuffles like the top to random shuffle, the transposing neighbors shuffle and the random transpositions shuffle, but also models for shuffles that people actually use for shuffling a deck of cards; the riffle shuffle and the overhand shuffle. The riffle shuffle (or the Gilbert–Shannon–Reeds shuffle) is a model for the most frequently used shuffling method. Aldous [1] showed that the cutoff

Received January 2005; revised July 2005.
[1]Supported in part by the Swedish Research Council.

for a deck of $n$ cards is given by $\frac{3}{2}\log_2 n$. Later Bayer and Diaconis [2] in a celebrated paper extended this to a very sharp result which, among other things, established that one needs about seven shuffles to mix an ordinary deck with 52 cards. The *overhand shuffle*, the topic of the present paper, is the shuffling technique where you gradually transfer the deck from, say, your right hand to you left hand by sliding off small packets from the top of the deck. Pemantle [5] introduced and analyzed a probabilistic model for the overhand shuffle and established upper and lower bounds of order $n^2 \log n$ and $n^2$, respectively, for the mixing time on a deck of $n$ cards. Here we will establish a lower bound of order $n^2 \log n$, thereby proving Pemantle's upper bound to be essentially tight. The main tool of our analysis is a technique introduced by Wilson [7] where one cleverly uses an eigenvalue close to 1 and a corresponding eigenvector of the transition matrix of the motion of a single card in order to find a useful eigenvector of the whole chain. This eigenvector exhibits a certain contraction property and the random variable obtained by evaluating it at the state at time $t$ can often be proved to be close to its expectation with high probability. In [7] this technique is used to prove an order $n^3 \log n$ lower bound for the transposing neighbors shuffle and in [6] a version for a lifted Markov chain is used to establish a lower bound of the same order for the Rudvalis shuffle. Mossel, Peres and Sinclair [4] use the technique to prove an $n \log n$ lower bound for the so called cyclic to random shuffle.

In Section 2 we review some basic definitions and describe Pemantle's model for the overhand shuffle and Wilson's technique. In Section 3 we state and prove the promised lower bound for the overhand shuffle. In order to do this, we will have to extend Wilson's technique slightly and use a vector that is not really but only very close to an eigenvector of the transition matrix. We hope that this extension turns out to be useful in other situations as well. In Section 4 we demonstrate how it can be used to give a simplified proof of Wilson's [6] $\Omega(n^3 \log n)$ lower bound for the Rudvalis shuffle.

## 2. Preliminaries.

2.1. *Basic definitions.* The standard way of measuring the distance between two probability measures $\mu$ and $\nu$ on a finite set $S$ is by the *total variation norm* defined by

$$\|\mu - \nu\| = \tfrac{1}{2}\sum_{s\in S}|\mu(S) - \nu(S)| = \max_{A \subseteq S}(\mu(A) - \nu(A)).$$

Let $\{X_t^n\}_{t=0}^{\infty}$, $n = 1, 2, 3, \ldots$, be a sequence of aperiodic irreducible Markov chains on state spaces $S^n$ and with stationary distributions $\pi^n$. Assume that $|S^n| \uparrow \infty$ in a natural way. A sequence $\{f(n)\}_{n=1}^{\infty}$ is said to be an *upper bound* for the mixing time of the sequence of Markov chains if

$$\lim_{n \to \infty} \|P(X_{f(n)}^n \in \cdot) - \pi^n\| = 0$$

and a *lower bound* if

$$\lim_{n \to \infty} \|P(X_{f(n)}^n \in \cdot) - \pi^m\| = 1.$$

When, for any $a > 0$, $(1-a)f(n)$ is a lower bound and $(1+a)f(n)$ is an upper bound, one says that $f(n)$ is a *cutoff* or a *threshold* for the mixing time. In card shuffling situations, as in the present paper, $S^n = S_n$, the set of permutations of $n$ cards, and usually $\pi^n$ is $U = U^n$, the uniform distribution on $S_n$.

2.2. *The overhand shuffle.* Pemantle's model for the overhand shuffle is parameterized by a probability $p \in (0, 1)$. The transition rule is the following: Each of the $n - 1$ slots between adjacent cards is, independently of the other slots, declared a *cutpoint* with probability $p$. This divides in a natural way the deck into subsequences, or packets, of cards. Reverse each of these packets without changing its position relative to the other packets. (This is not what really happens when doing an overhand shuffle; what we get is rather the reversed deck compared to what one would have got from the overhand shuffle. This of course does not matter for the mixing time.)

For this Markov chain, Pemantle obtained an $O(n^2 \log n)$ upper bound by constructing a coupling for which he was able to demonstrate that one can find a constant $C < \infty$ such that each card is with high probability matched within $Cn^2 \log n$ steps. Intuitively, one can understand this by observing that a particular card makes a random walk which is essentially symmetric [at least when the card is not at $O(1)$ distance from the ends of the deck] with step size variance of constant order. Therefore, one should be able to match the typical card within $O(n^2)$ steps and all cards within $O(n^2 \log n)$ steps. Studying the random walks of single cards also quickly yields an $\Omega(n^2)$ lower bound since before this time too many of the cards are too close to their starting positions.

2.3. *Wilson's technique.* The key to the technique is the following lemma from [7]. Later on, an extension of the lemma will be needed together with corresponding adjustments of its proof. We prove the original lemma here in order to make the later proof of the extension easier to digest; doing both at once may hide the essential ideas behind technical details.

LEMMA 2.1 (The lower bound lemma). *Let $\{X_t\}_{t=0}^{\infty}$ be an aperiodic irreducible Markov chain with state space $S$ and stationary distribution $\pi$.*

*Assume that* $\Phi\colon S \to \mathbb{R}$ *and* $\gamma \in (0, \frac{1}{2})$ *are such that*

$$\mathbb{E}[\Phi(X_{t+1})|X_t] = (1-\gamma)\Phi(X_t).$$

*(I.e., $1-\gamma$ is an eigenvalue for the transition matrix and $\Phi$ is a corresponding eigenvector.) Suppose that $R \geq \max_{s \in S} \mathbb{E}[(\Phi(X_{t+1}) - \Phi(X_t))^2 | X_t = s]$, put $\hat{\Phi} = \max_{s \in S} \Phi(s)$ and let the chain start from a state $s_0$ such that $\Phi(s_0) = \hat{\Phi}$.*

*Fix $\varepsilon > 0$ and put*

$$T = \frac{\log \hat{\Phi} - (1/2)\log(4R/(\gamma\varepsilon))}{-\log(1-\gamma)}.$$

*Then for $t \leq T$,*

$$\|P(X_t \in \cdot) - \pi\| \geq 1 - \varepsilon.$$

The lemma may be adopted to situations with complex-valued eigenvalues and eigenvectors as well, but we will not need that here; the overhand shuffle clearly defines a reversible Markov chain.

PROOF OF LEMMA 3.2. By induction, $\mathbb{E}\Phi(X_t) = (1-\gamma)^t \Phi(s_0)$ and, consequently, $\mathbb{E}\Phi(X_\infty) = 0$. Put $\Delta\Phi = \Phi(X_{t+1}) - \Phi(X_t)$. Then

$$\mathbb{E}[\Phi(X_{t+1})^2 | X_t] = (1-2\gamma)\Phi(X_t)^2 + \mathbb{E}[(\Delta\Phi)^2 | X_t] \leq (1-2\gamma)\Phi(X_t)^2 + R.$$

Using induction and that $\gamma \leq 1/2$, we get

$$\mathbb{E}[\Phi(X_t)^2] \leq (1-2\gamma)^t \Phi(s_0)^2 + \frac{R}{2\gamma}.$$

Thus, since $\gamma \leq 1/2$,

$$\operatorname{Var}\Phi(X_t) \leq \Phi(s_0)^2((1-2\gamma)^t - (1-\gamma)^{2t}) + \frac{R}{2\gamma} \leq \frac{R}{2\gamma}.$$

By Chebyshev's inequality,

$$P\left(|\Phi(X_t) - \mathbb{E}\Phi(X_t)| \geq \sqrt{\frac{R}{\gamma\varepsilon}}\right) \leq \frac{\varepsilon}{2}$$

and by letting $t \to \infty$, it also follows that

$$P\left(|\Phi(X_\infty)| \geq \sqrt{\frac{R}{\gamma\varepsilon}}\right) \leq \frac{\varepsilon}{2}.$$

Now if $t$ is such that $\mathbb{E}\Phi(X_t) \geq 2\sqrt{R/(\gamma\varepsilon)}$, then we get that $P(\Phi(X_t) \leq \sqrt{R/(\gamma\varepsilon)}) \leq \varepsilon/2$ and so

$$\|P(X_t \in \cdot) - P(X_\infty \in \cdot)\| \geq 1 - \varepsilon,$$

as desired. It remains to determine for which $t$ we have $\mathbb{E}\Phi(X_t) = (1-\gamma)^t \hat{\Phi} \geq 2\sqrt{R/(\gamma\varepsilon)}$. Some simple algebraic manipulations show that this is the case when $t \leq T$. □

There is more to using the lower bound lemma than just the lemma itself: In order to derive good lower bounds, one needs to find an eigenvalue sufficiently close to 1 and a useful corresponding eigenvector. It is not at all obvious how to do that. In [7], and also in [6] and [4], one may find some useful hints on how this may be done.

**3. Lower bound for the overhand shuffle.** Let $X_t$ denote the state of the deck at time $t$ and let $Z_t = Z_t^i$ denote the position of card $i$ at time $t$. For simplicity, we shall do the the case $p = 1/2$ and leave the straightforward generalization to the reader. To further simplify our calculations, we will to begin with assuming a *circular deck convention*, that is, regarding the top card and the bottom card as being next to each other. This means that the top packet and the bottom packet may be treated as a single packet (if the slot between that top and the bottom card happens to be a cutpoint). Then any single card performs a symmetric random walk where the probability that the card moves $k$ steps to the right (modulo $n$) at a single shuffle is $\frac{1}{3}(\frac{1}{2})^{|k|}$, $k \in \mathbb{Z}$. (A negative $k$, of course, means that the card moves to the left. The given step size distributions, not entirely correct because of the remote possibility that, at a given shuffle, there are no cutpoints at all or only one cutpoint. However, the probability that this ever happens in the time scale we are interested in is so small that we do not need to care about this. If we want to make the formula completely correct, we may just redefine the shuffle on this event in a suitable way.) We start with the circular deck convention since it produces a much cleaner proof than the "real" overhand shuffle. At the end of this section the circular convention will be removed via a slight extension of the lower bound lemma.

To find an eigenvector for the $\{Z_t\}$-chain, we use the cosine function [note that $\cos(2\pi k/n)$ is well defined on $\mathbb{Z}_n$]:

$$\mathbb{E}\left[\cos\frac{2\pi Z_{t+1}}{n}\Big|X_t\right]$$

$$= \frac{1}{3}\cos\frac{2\pi Z_t}{n} + \frac{1}{3}\sum_{k=1}^{\infty}\frac{1}{2^k}\left(\cos\frac{2\pi(Z_t-k)}{n} + \cos\frac{2\pi(Z_t+k)}{n}\right)$$

$$= \frac{1}{3}\left(\cos\frac{2\pi Z_t}{n}\right)\left(1 + 2\sum_{k=1}^{\infty}\frac{\cos(2\pi k/n)}{2^k}\right) = (1-\gamma)\cos\frac{2\pi Z_t}{n},$$

where $\gamma = (1+o(1))8\pi^2 n^{-2}$, which follows from the fact that, for $k = o(n)$, one has $\cos(2\pi k/n) = 1 - 2\pi^2 k^2 n^{-2} + O(k^4 n^{-4})$, and that $\sum_{k=1}^{\infty} k^2 2^{-k} = 6$.

The second equality follows from the standard trigonometric formula for the cosine of a sum.

Now put $m = \lfloor n/2 \rfloor$ and

$$\Phi(X_t) = \sum_{i=1}^{m} \cos \frac{2\pi Z_t^i}{n}.$$

Then $\Phi$ is an eigenvector corresponding to the eigenvalue $1 - \gamma$ for the transition matrix of the whole deck, and $\hat{\Phi} = Cn$ for $C = \Theta(1)$. In order to use the lower bound lemma, we need to bound $\mathbb{E}[(\Phi(X_{t+1}) - \Phi(X_t))^2 | X_t]$. Write

$$Y_i = \cos \frac{2\pi Z_{t+1}^i}{n} - \cos \frac{2\pi Z_t^i}{n}.$$

Then

$$\mathbb{E}[(\Phi(X_{t+1}) - \Phi(X_t))^2 | X_t] = \mathbb{E}\left[\left(\sum_{i=1}^{m} Y_i\right)^2 \Big| X_t\right] = \sum_{i=1}^{m} \sum_{j=1}^{m} \mathbb{E}[Y_i Y_j | X_t].$$

Let $A$ denote the event that no two cutpoints of the $(t+1)$st shuffle are more than $10 \log n$ steps apart. It is easily verified that $P(A) = 1 - o(n^{-5})$.

Consider a given pair $(i, j)$ of cards that at time $t$ are more than $20 \log n$ apart. Ideally, if card $i$ and card $j$ were infinitely far apart at time $t$, then $Y_i$ and $Y_j$ would have been independent. Put $Y_i'$ and $Y_j'$ for such idealized independent random variables. By the definition of $A$, we can couple $Y_i$ and $Y_j$ with $Y_i'$ and $Y_j'$ in such a way that $Y_i = Y_i'$ and $Y_j = Y_j'$ on $A$. (Such a coupling can be made by coupling the cutpoints for the two cases in a natural way.) Thus,

$$\mathbb{E}[Y_i Y_j | X_t] = \mathbb{E}[Y_i' Y_j' I_A | X_t] + o(n^{-5})$$
$$= \mathbb{E}[Y_i' | X_t] \mathbb{E}[Y_j' | X_t] + o(n^{-5})$$
$$= \mathbb{E}[Y_i | X_t] \mathbb{E}[Y_j | X_t] + o(n^{-5}).$$

By mimicking the above calculations, it follows that $\mathbb{E}[Y_i | X_t] = O(n^{-2})$ and so

$$\mathbb{E}[Y_i Y_j | X_t] = O(n^{-4}).$$

Using this, it follows that

$$\mathbb{E}[(\Phi(X_{t+1}) - \Phi(X_t))^2 | X_t]$$
$$= \sum_{i=1}^{m} \sum_{j : |Z_t^i - Z_t^j| \leq 20 \log n} \mathbb{E}[Y_i Y_j | X_t] + O(n^{-2})$$
$$= \sum_{i=1}^{m} \sum_{j : |Z_t^i - Z_t^j| \leq 20 \log n} \mathrm{Cov}(Y_i, Y_j | X_t) + O(n^{-2}) \leq 20 V n \log n,$$

where
$$V = \max_{s \in S_n} \text{Var}(Y_i|X_t = s) \leq \max_{s \in S_n} \mathbb{E}[Y_i^2|X_t = s]$$
$$= \frac{8}{3} \sum_{k=1}^{\infty} \frac{\sin^2 2\pi k/n}{2^k} \leq 100 n^{-2}$$

for large enough $n$. Hence,
$$\mathbb{E}[(\Phi(X_{t+1}) - \Phi(X_t))^2 | X_t] \leq 2000 n^{-1} \log n$$

and we may apply the lower bound lemma with $R = 2000 n^{-1} \log n$. In order to have $\varepsilon \to 0$, let $\varepsilon = (\log n)^{-1}$. The lower bound lemma now yields the lower bound

$$T = (1 + o(1)) \frac{n^2}{8\pi^2} \left( \log(Cn) - \frac{1}{2} \log \frac{8000 n (\log n)^2}{8\pi^2} \right)$$
$$= (1 + o(1)) \frac{1}{16\pi^2} n^2 \log n.$$

Finally, for general $p$, the analogous calculations using that $\sum_{k=1}^{\infty} j^2 (1-p)^j = (2-p)(1-p)/p^3$ give the following:

THEOREM 3.1. *A lower bound for the overhand shuffle with circular deck convention and with parameter $p \in (0,1)$ is given by*

$$\frac{p^2(2-p)}{8\pi^2(1-p^2)} n^2 \log n.$$

Next we remove the circular deck convention. We will again concentrate on the case $p = 1/2$. The vector $\Phi$ above is now no longer really an eigenvector for $\{X_t\}$, but only very close to one; in particular, the cards near the ends of the deck do not behave well. However, the feeling is that the small error we make if we still assume that $\Phi$ is an eigenvector should not be large enough to upset things too much. This is made precise by the following extension of the lower bound lemma:

LEMMA 3.2 (The extended lower bound lemma). *Let the setting be as in the lower bound lemma with the exception that $\Phi$ is not an eigenvector for the transition matrix, but there is a positive number $\rho$ such that*

$$(1-\gamma)\Phi(X_t) - \rho \leq \mathbb{E}[\Phi(X_{t+1})|X_t] \leq (1-\gamma)\Phi(X_t) + \rho$$

*a.s. Then, for $t \leq T$, $\|P(X_T \in \cdot) - \pi\| \geq 1 - \varepsilon$, where*

$$T = \frac{\log \hat{\Phi} - \log(2\rho/\gamma + \sqrt{4(R + 6\rho\hat{\Phi})/(\gamma\varepsilon)})}{-\log(1-\gamma)}.$$

PROOF. The proof is basically a repetition of the proof of the lower bound lemma, with the proper adjustments: By induction, we get

$$(1-\gamma)^t \Phi(s_0) - \frac{\rho}{\gamma} \leq \mathbb{E}\Phi(X_t) \leq (1-\gamma)^t \Phi(s_0) + \frac{\rho}{\gamma},$$

that is, for some $\alpha \in [-1, 1]$,

$$\mathbb{E}\Phi(X_t) = (1-\gamma)^t \Phi(s_0) - \frac{\alpha\rho}{\gamma}$$

and by letting $t \to \infty$, we also get

$$-\frac{\rho}{\gamma} \leq \mathbb{E}\Phi(X_\infty) \leq \frac{\rho}{\gamma}.$$

We have

$$\mathbb{E}[\Phi(X_t)^2 | X_t] \leq \Phi(X_t)^2 + 2\Phi(X_t)\mathbb{E}[\Delta\Phi | X_t] + R$$
$$\leq (1-2\gamma)\Phi(X_t)^2 + 2\rho|\Phi(X_t)| + R$$
$$\leq (1-2\gamma)\Phi(X_t)^2 + 2\rho\hat{\Phi} + R.$$

By induction,

$$\mathbb{E}[\Phi(X_t)^2] \leq (1-2\gamma)^t \Phi(s_0)^2 + \frac{R + 2\rho\hat{\Phi}}{2\gamma}.$$

Thus,

$$\operatorname{Var}\Phi(X_t) \leq (1-2\gamma)^t \Phi(s_0)^2 + \frac{R + 2\rho\hat{\Phi}}{2\gamma} - \left((1-\gamma)^t \Phi(s_0) - \frac{\alpha\rho}{\gamma}\right)^2$$
$$\leq \frac{R + 2\rho\hat{\Phi}}{2\gamma} + \frac{2|\alpha|\rho\hat{\Phi}(1-\gamma)^t}{\gamma} \leq \frac{R + 6\rho\hat{\Phi}}{2\gamma}.$$

By Chebyshev's inequality,

$$P\left(|\Phi(X_t) - \mathbb{E}\Phi(X_t)| \geq \sqrt{\frac{R + 6\rho\hat{\Phi}}{\gamma\varepsilon}}\right) \leq \frac{\varepsilon}{2}$$

and by letting $t \to \infty$, using that $\mathbb{E}\Phi(X_\infty) \leq \rho/\gamma$,

$$P\left(|\Phi(X_\infty)| \geq \frac{\rho}{\gamma} + \sqrt{\frac{R + 6\rho\hat{\Phi}}{\gamma\varepsilon}}\right) \leq \frac{\varepsilon}{2}.$$

Thus, $\|P(X_t \in \cdot) - \pi\| \geq 1 - \varepsilon$ when $t$ is such that $\mathbb{E}\Phi(X_t) \geq \rho/\gamma + 2\sqrt{(R + 6\rho\hat{\Phi})/(\gamma\varepsilon)}$. Since $\mathbb{E}\Phi(X_t) \geq (1-\gamma)^t \hat{\Phi} - \rho/\gamma$, this holds when $t \leq T$.
□

In our case $\gamma = \Theta(n^{-2})$, $R = O(n^{-1}\log n)$ (the above arguments apply equally well to the present setting without the circular deck convention) and $\hat{\Phi} = O(n)$. Hence, once it has been shown that $\rho = O(n^{-2+\delta})$ for any $\delta > 0$, then we get the same lower bound that we got with the circular deck convention. In fact, with some care it can be shown that $\rho = O(n^{-2})$, but things get considerably easier if we settle for proving that $\rho = O(n^{-2}(\log n)^3)$.

As above, the probability that in a given shuffle the position of the first cutpoint is larger than $10\log n$ or that the position of the last cutpoint is less than $n - 10\log n$ is $o(n^{-5})$. Therefore, when $X_t$ is such that a given card is in a position $Z_t = k$ with $10\log n < k < n - 10\log n$, the movement of that card can be coupled with the movement of a card under the circular deck convention so that the two movements coincide with probability $1 - o(n^{-5})$. Thus,

$$E\left[\cos\frac{2\pi Z_{t+1}}{n}\Big|X_t\right] = (1-\gamma)\cos\frac{2\pi Z_t}{n} + o(n^{-5}).$$

Now assume that $k \leq 10\log n$ (the case $k \geq n - 10\log n$ is treated in the same way). Then

$$\cos\frac{2\pi k}{n} = 1 - O\left(\frac{(\log n)^2}{n^2}\right).$$

Since $P(|Z_{t+1} - Z_t| \geq 2\log_2 n) \leq n^{-2}$, we get

$$\mathbb{E}\left[\cos\frac{2\pi Z_{t+1}}{n}\Big|X_t\right] \geq \left(1 - \frac{1}{n^2}\right)\cos\frac{2\pi(k + 2\log_2 n)}{n} - \frac{1}{n^2}$$

$$= 1 - O\left(\frac{(\log n)^2}{n^2}\right)$$

$$= (1-\gamma)\cos\frac{2\pi Z_t}{n} + O\left(\frac{(\log n)^2}{n^2}\right).$$

Since

$$\Phi(X_t) = \sum_{i=1}^{\lfloor n/2 \rfloor} \cos\frac{2\pi Z_t^i}{n},$$

it follows on summation that

$$\mathbb{E}[\Phi(X_{t+1})|X_t] = (1-\gamma)\Phi(X_t) + 20\log n \cdot O\left(\frac{(\log n)^2}{n^2}\right) + \frac{n}{2}o(n^{-5})$$

$$= (1-\gamma)\Phi(X_t) + O\left(\frac{(\log n)^3}{n^2}\right)$$

and so $\rho = O((\log n)^3 n^{-2})$ as desired.

Again, the case with a general $p \in (0,1)$ is completely analogous and we arrive at the following:

THEOREM 3.3. *A lower bound for the overhand shuffle with parameter $p \in (0,1)$ is given by*

$$\frac{p^2(2-p)}{8\pi^2(1-p^2)} n^2 \log n.$$

**4. The Rudvalis shuffle.** The (inverse) Rudvalis shuffle is described by the following transition rule: With probability $1/2$, move the bottom card to the top of the deck and with probability $1/2$, move the card next to the bottom card to the top of the deck. The model was proposed by Arunas Rudvalis as a shuffle with very slow convergence to uniformity. In [3] it was shown that the Rudvalis shuffle mixes in $O(n^3 \log n)$ steps. For a long time the best known lower bound was $\Omega(n^3)$, until Wilson [6] recently found an $\Omega(n^3 \log n)$ bound and thereby established the mixing time as $\Theta(n^3 \log n)$.

Applying the lower bound lemma directly to the Rudvalis shuffle is difficult, the reason being that the nontrivial eigenvalues of the chain described by the motion of a single card are complex and deviate from 1 by $\Omega(n^{-1})$. Wilson gets around this problem by lifting the Rudvalis chain to a larger state space by letting time be a part of this larger state space. He then extends the lower bound lemma in a way that lets him make conclusions about the mixing time of the Rudvalis chain from the behavior of the lifted chain.

Here we shall reprove Wilson's lower bound by applying our extended lower bound lemma above to the $n$-step transition matrix of the Rudvalis shuffle. More precisely, let $X_t$ as before denote the state of the deck at time $t$ and $Z_t^i$ the position of card $i$ at time $t$. Then with $\gamma = (1+o(1))4\pi^2 n^{-2}$ and $\Phi(X_t) = \sum_{i=1}^{\lfloor n/2 \rfloor} \cos(2\pi Z_t^i/n)$, the pair $(1-\gamma, \Phi(X_t))$ is sufficiently close to an eigenvalue/eigenvector pair for the chain $\{X_{ns}\}_{s=1}^\infty$ for the extended lower bound lemma to be invoked. Due to the fact that we are content with finding a vector which is not really but only close enough to an eigenvector, we are able to avoid some technical difficulties. For example, we will not need to lift the chain and we will not have to deal with any complex algebra whatsoever.

Let $p_{k,j}$ denote the probability that a card starting from position $k$ finds itself in position $j$ after $n$ steps of the Rudvalis shuffle. If $k \notin \{n-1, n\}$, then the first $n-k-1$ steps deterministically take the card to position $n-1$. From there it takes the card a geometric$(1/2)$ number $G$, say, of steps for the card to move to position 1. After that, provided it happens before all the $n$ steps have been taken, the remaining $k+1-G$ steps deterministically shift the card to position $k+2-G$. We get

$$p_{k,n} = (\tfrac{1}{2})^{k+1}$$

and for $j = 1, 2, \ldots, k+1$,
$$p_{k,j} = (\tfrac{1}{2})^{k+2-j}.$$
A special treatment along the same lines of the case $k \in \{n-1, n\}$ gives
$$p_{k,n} = p_{k,1} = \tfrac{1}{4} + (\tfrac{1}{2})^n$$
and for $j = 2, 3, \ldots, n-1$,
$$p_{k,j} = (\tfrac{1}{2})^{n+1-j}.$$
The analysis is now very similar to that of the overhand shuffle above. Suppose $X_{ns}$ is such that $Z_{nt} = Z_{nt}^i = k$ for a given card $i$. When $k \in \{n, n-1\}$ or $k \leq 10 \log n$, $\cos(2\pi k/n) = 1 - O((\log n)^2 n^{-2})$. Since $P(|Z_{n(s+1)} - Z_{ns}| \geq 2 \log_2 n) \leq n^{-2}$,
$$\mathbb{E}\left[\cos \frac{2\pi Z_{n(s+1)}}{n} \Big| X_{ns}\right] = 1 - O\left(\frac{(\log n)^2}{n^2}\right)$$
$$= (1-\gamma) \cos \frac{2\pi Z_{ns}}{n} + O\left(\frac{(\log n)^2}{n^2}\right).$$
Now assume that $k \in \{\lceil 10 \log n \rceil, \ldots, n-2\}$. Then
$$\mathbb{E}\left[\cos \frac{2\pi Z_{n(s+1)}}{n} \Big| X_{ns}\right]$$
$$= \sum_{j=-k}^{1} \frac{1}{2^{2-j}} \cos \frac{2\pi(k+j)}{n} + \frac{1}{2^{k+2}}$$
$$= \left(\sum_{j=-k}^{1} \frac{1}{2^{2-j}} \cos \frac{2\pi j}{n}\right) \cos \frac{2\pi k}{n} - \left(\sum_{j=-k}^{1} \frac{1}{2^{2-j}} \sin \frac{2\pi j}{n}\right) \sin \frac{2\pi k}{n}$$
$$+ o(n^{-5}).$$
The first of the three terms equals $(1-\gamma) \cos(2\pi k/n)$. For the second term,
$$-\sum_{j=-k}^{1} \frac{1}{2^{2-j}} \sin \frac{2\pi j}{n} = \sum_{j=-1}^{k} \frac{1}{2^{2+j}} \sin \frac{2\pi j}{n}$$
$$= \sum_{j=-1}^{k} \frac{1}{2^{2+j}} \left(\frac{2\pi j}{n} + O(n^{-3})\right)$$
$$= \frac{2\pi}{n} \sum_{j=-1}^{\infty} \frac{j}{2^{2+j}} + O(n^{-3}) = O(n^{-3}),$$
where the last equality follows from the fact that the last sum is 0, which can, for example, be seen by noting that the sum coincides with $-2$ plus the expectation of a geometric(1/2) random variable. On summing the three terms and then over the cards $1, 2, \ldots, \lfloor n/2 \rfloor$, we get

$$\mathbb{E}[\Phi(X_{n(s+1)})|X_{ns}] = (1-\gamma)\Phi(X_{ns}) + O\left(\frac{(\log n)^3}{n^2}\right).$$

In order to apply the extended lower bound lemma, we need to bound $\mathbb{E}[(\Phi(X_{n(s+1)}) - \Phi(X_{ns}))^2|X_{ns}]$. Put

$$U_t^i = \cos\frac{2\pi(Z_t^i - t)}{n}.$$

Then, with $m = \lfloor n/2 \rfloor$,

$$\Phi(X_{n(s+1)}) - \Phi(X_{ns}) = \sum_{i=1}^{m}(U_{n(s+1)}^i - U_{ns}^i)$$

$$= \sum_{i=1}^{m}\sum_{t=ns+1}^{n(s+1)}(U_t^i - U_{t-1}^i)$$

$$= \sum_{t=ns+1}^{n(s+1)}\sum_{i=1}^{m}(U_t^i - U_{t-1}^i).$$

Since the terms for different $t$'s are concerned with different shuffles, these are independent. Also, since all cards but at most two are at each shuffle shifted one step down the deck, at most two terms of the inner sum can change, and if so, the change is bounded by $2/n$. Taking these facts into account, we get

$$\mathrm{Var}(\Phi(X_{n(s+1)}) - \Phi(X_{ns})|X_{ns}) \leq m \cdot \frac{4}{n^2} \leq \frac{2}{n}.$$

Adding this to

$$(\mathbb{E}[\Phi(X_{n(s+1)}) - \Phi(X_{ns})|X_t])^2 = \left(\gamma\Phi(X_{ns}) + O\left(\frac{(\log n)^3}{n^2}\right)\right)^2 = O(n^{-2}),$$

we get

$$\mathbb{E}[(\Phi(X_{n(s+1)}) - \Phi(X_{ns}))^2|X_{ns}] = O(n^{-1}).$$

We may thus apply the extended lower bound lemma with $\gamma = (1+o(1))4\pi^2 n^{-2}$, $R = O(n^{-1})$ and $\rho = O((\log n)^3 n^{-2})$. Doing so, we arrive at the lower bound $\frac{1}{8\pi^2}n^2\log n$ for $\{X_{ns}\}_{s=1}^{\infty}$. For the Rudvalis shuffle, this gives the lower bound

$$\frac{1}{8\pi^2}n^3\log n,$$

as desired.

**Acknowledgment.** I am grateful to the referee for the thorough reading of the manuscript.

## REFERENCES

[1] ALDOUS, D. (1983). Random walks on finite groups and rapidly mixing Markov chains. *Séminaire de Probabilités XVII. Lecture Notes in Math.* **986** 243–297. Springer, Berlin. MR770418
[2] BAYER, D. and DIACONIS, P. (1992). Trailing the dovetail shuffle to its lair. *Ann. Appl. Probab.* **2** 294–313. MR1161056
[3] HILDEBRAND, M. (1990). Rates of convergence of some random processes on finite groups. Ph.D. thesis, Harvard Univ.
[4] MOSSEL, E., PERES, Y. and SINCLAIR, A. (2004). Shuffling by semi-random transpositions. *FOCS 2004.*
[5] PEMANTLE, R. (1989). Randomization time for the overhand shuffle. *J. Theoret. Probab.* **2** 37–49. MR981762
[6] WILSON, D. B. (2003). Mixing time of the Rudvalis shuffle. *Electron. Comm. Probab.* **8** 77–85. MR1987096
[7] WILSON, D. B. (2004). Mixing times of lozenge tiling and card shuffling Markov chains. *Ann. Appl. Probab.* **14** 274–325. MR2023023

DEPARTMENT OF MATHEMATICS
CHALMERS UNIVERSITY OF TECHNOLOGY
S-412 96 GÖTEBORG
SWEDEN
E-MAIL: jonasson@math.chalmers.se